\documentclass[a4paper,12pt]{article}
\usepackage{comment}
\usepackage{cite}
\usepackage{amsmath}
\usepackage{amssymb}
\usepackage{amsfonts}
\usepackage[T1]{fontenc}
\usepackage[utf8]{inputenc}
\usepackage{graphicx}
\usepackage{fancyhdr}
\usepackage{float}
\usepackage{xcolor}
\usepackage{authblk}
\usepackage{mathrsfs}
\usepackage{empheq}
\usepackage[hyphens]{url}
\usepackage{hyperref} 
\usepackage[]{breakurl}
\usepackage{amsmath}
\usepackage{amssymb}

\usepackage{geometry}
\begin{document}

\title{Bounds on the $p$-adic valuation of the factorial, hyperfactorial and superfactorial}

\author[$\dagger$]{Jean-Christophe {\sc Pain}$^{1,2,}$\footnote{jean-christophe.pain@cea.fr}\\
\small
$^1$CEA, DAM, DIF, F-91297 Arpajon, France\\
$^2$Universit\'e Paris-Saclay, CEA, Laboratoire Mati\`ere en Conditions Extr\^emes,\\ 
F-91680 Bruy\`eres-le-Ch\^atel, France
}

\date{}

\maketitle

\begin{abstract}
In this article, we investigate the $p$-adic valuation $\nu_p$ of quantities such as the factorial $n!$, the hyperfactorial $H(n)$ or the superfactorial $\mathrm{sf}(n)$. In particular, we obtain simple bounds (both upper and lower) for $\nu_p$, using the Legendre-de Polignac formula. Other iterated quantities such as the Berezin function, are also considered. Beyond  their ``recreational'' character, such quantities, often related to very large numbers, may find applications for cryptography purposes. Finally, lower and upper bounds for the $p$-adic valuation of Stirling numbers of the first kind and Catalan numbers are briefly discussed.
\end{abstract}

\section{Introduction}

Besides usual quantities such as $n!$ or the double factorial $n!!=n(n-2)(n-4)\cdots$, some more complex quantities were introduced in the framework of the study of very large numbers, such as the hyperfactorial:
\begin{equation}
    H(n)=\prod_{k=1}^nk^k
\end{equation}
or the superfactorial:
\begin{equation}
    \mathrm{sf}(n)=\prod_{k=1}^nk!.
\end{equation}
The $p$-adic valuation (or $p$-adic order) of an integer $n$, denoted $\nu_{p}(n)$ \cite{Niven1991}, is the exponent of the highest power of the prime number $p$ that divides $n$. Equivalently, $\nu_{p}(n)$ is the exponent to which $p$ appears in the prime factorization of $n$. The $p$-adic valuation of an integer $n\in\mathbb{N}$ is defined to be
\begin{equation}
    \nu_{p}(n)=\mathrm{max} \{k\in \mathbb{N}:p^{k}\mid n\}.
\end{equation}
In additive and algebraic number theory, the Skolem-Mahler-Lech theorem \cite{Skolem1933,Mahler1935,Lech1953} states that if a sequence of numbers satisfies a linear difference equation, then with finitely many exceptions the positions at which the sequence is zero form a regularly repeating pattern\footnote{This result is named after Skolem (who proved the theorem for sequences of rational numbers), Mahler (who proved it for sequences of algebraic numbers), and Lech (who proved it for sequences whose elements belong to any field of characteristic 0).} Its known proofs use $p$-adic analysis and are non-constructive. 

Another application of $p$-adic numbers is the proof, given in Serre's course in arithmetic \cite{Serre1973}, that a natural number is expressible as the sum of squares if and only if it is not of the form $4a(8b+11)$ for some $a,b\in\mathbb{N}^2$.
 
In section \ref{sec2}, we discuss simple bounds of $\nu_p$ for the $p$-adic valuation of the factorial, hyperfactorial and superfactorial. Our bounds are mainly based on the Legendre-de Polignac formula
\begin{equation}
    n!=\prod_{i=1}^{\pi(n)}p_i^{\sum_{j=1}^{\lfloor\log_{p_i} n\rfloor}\Bigl\lfloor\frac{n}{p_i^j}\Bigl\rfloor},
\end{equation}
where $\pi(n)$ is the prime counting function, i.e., the number of prime numbers less than or equal to $n$.

Other iterated quantities such as the Berezin function \cite{Berezin1987}:
\begin{equation}
    n\$=n!^{n!^{n!^{\cdot^{\cdot^{\cdot}}}}},
\end{equation}
$n!$ appearing $n!$ times, and the functions

\begin{equation}
    F_1(n)=\prod_{k=1}^nk^{k!},\;\;\;\; F_2(n)=\prod_{k=1}^nk!^{k}\;\;\;\;\mathrm{and}\;\;\;\;F_3(n)=\prod_{k=1}^nk!^{k!}
\end{equation}
are evoked in section \ref{sec3} and additional relations for the $p$-adic valuation are given in Appendix \ref{appA}, still in the scope of obtaining inequalities.

\section{Applying the Legendre-de Polignac formula}\label{sec2}

\subsection{Factorial}

Let us consider $n\in\mathbb{N}^*$. We have, for $k\in\mathbb{N}^*$, according to the Legendre-de Polignac formula \cite{Legendre}:
\begin{equation}
    \nu_p(n!)=\sum_{k=1}^{\infty}\Bigl\lfloor\frac{n}{p^k}\Bigl\rfloor
\end{equation}
and since
\begin{equation}
    \Bigl\lfloor\frac{n}{p^k}\Bigl\rfloor\leq\frac{n}{p^k} 
\end{equation}
we get
\begin{equation}
    \nu_p(n!)\leq\sum_{k=1}^{\infty}\frac{n}{p^k}.
\end{equation}
We have also
\begin{equation}
    \sum_{k=1}^{\infty}\frac{n}{p^k}=\frac{n}{p}\sum_{k=0}^{\infty}\frac{1}{p^k}=\frac{n}{p}\frac{1}{(1-1/p)}=\frac{n}{p-1}.
\end{equation}
Thus one has
\begin{equation}\label{majo1}
    \nu_p(n!)\leq\frac{n}{p-1}.
\end{equation}
However, a better upper bound can be obtained, together with a lower bound \cite{Marques2012}. Indeed, one has \cite{Legendre}:
\begin{equation}\label{d2}
    \nu_p(n!)=\frac{n-s_p(n)}{p-1},
\end{equation}
where $s_p(n)$ is the sum of digits of $n$ in base $p$. Since $n$ has $\lfloor\log n/\log p\rfloor$ digits in base $p$, and each digit is at most $p-1$, we get
\begin{equation}\label{ine1}
    1\leq s_p(n)\leq(p-1)\left(\Bigl\lfloor\frac{\log n}{\log p}\Bigl\rfloor +1\right).
\end{equation}
Combining (\ref{ine1}) with expression (\ref{d2}) gives
\begin{equation}
    \frac{n}{p-1}-\Bigl\lfloor\frac{\log n}{\log p}\Bigl\rfloor-1\leq\nu_p(n!)\leq\frac{n-1}{p-1}
\end{equation}
or equivalently
\begin{equation}\label{mainineq}
        \frac{n}{p-1}-\Bigl\lfloor\log_p n\Bigl\rfloor-1\leq\nu_p(n!)\leq\frac{n-1}{p-1}.
\end{equation}

\subsection{Hyperfactorial}

Lots of studies have been done about the hyperfactorial function. In particular, Glaisher and Kinkelin found its asymptotic behaviour as $n$ approaches infinity. This led to the introduction of a constant, named after them \cite{Kinkelin1860,Glaisher1878,Glaisher1894,Finch2003,Chen2012}, which has a lot of expressions using the Euler Gamma function and the Riemann Zeta function.

The hyperfactorial of a positive integer $n$ is the product of the numbers $1^{1}, 2^{2}, \dots , n^{n}$, that is,
\begin{equation}\label{hyperfact}
    H(n)=1^{1}\cdot 2^{2} \cdots n^{n}=\prod_{i=1}^{n}i^{i}=n^{n}H(n-1).
\end{equation}
Following the usual convention for the empty product, the hyperfactorial of 0 is 1 ($H(0)=1$). The first values are 1, 1, 4, 108, 27648, 86400000, 4031078400000, 3319766398771200000, \emph{etc.}. It is worth mentioning that the hyperfactorial can be generalized to complex numbers, with
\begin{equation}
     H(z-1)\,G(z)=e^{(z-1)\log\Gamma(z)},
\end{equation}
where G is the Barnes G-function. The hyperfactorial has also the integral representation \cite{Mathematica}:
\begin{equation}
    H(z)=\frac{1}{(2\pi)^{z/2}}\exp\left[\binom{z+1}{2}+\int_0^z\log(t!)\,\mathrm{d}t\right] 	
\end{equation}
and the closed-form expression
\begin{equation}
    H(z-1)=\exp\left[\zeta'(-1,z+1)-\zeta'(-1)\right] 	
\end{equation}
for $\Re(z)>0$, where $\zeta(z)$ is the Riemann zeta function, $\zeta'(z)$ its derivative and $\zeta'(a,z)$ is the derivative of the Hurwitz zeta function with respect to the first argument.

In case of the hyperfactorial (see Eq. (\ref{hyperfact})), we have \cite{Onnis}:
\begin{equation}\label{basis}
    \nu_p\left[H(n)\right]=p\Bigl\lfloor\frac{n}{p}\Bigl\rfloor\,\nu_p(n!)-p\sum_{k=1}^{\Bigl\lfloor\frac{n}{p}\Bigl\rfloor-1}\nu_p\left[(p\times k)!\right],
\end{equation}
which can be proven by induction (see Ref. \cite{Onnis}). Using the Legendre-de Polignac formula, one has also
\begin{equation}
    \nu_p\left[H(n)\right]=p\Bigl\lfloor\frac{n}{p}\Bigl\rfloor\sum_{k=1}^{\infty}\Bigl\lfloor\frac{n}{p^k}\Bigl\rfloor-p\sum_{k=1}^{\Bigl\lfloor\frac{n}{p}\Bigl\rfloor-1}\sum_{l=1}^{\infty}\Bigl\lfloor\frac{i}{p^{k-1}}\Bigl\rfloor.
\end{equation}
Combining Eq. (\ref{basis}) with (\ref{mainineq}), one obtains
\begin{equation}
    \nu_p\left[H(n)\right]\geq p\Bigl\lfloor\frac{n}{p}\Bigl\rfloor\left(\frac{n}{p-1}-\Bigl\lfloor\log_p(n)\Bigl\rfloor-1\right)-p\sum_{k=1}^{\Big\lfloor\frac{n}{p}\Bigl\rfloor-1}\frac{pk-1}{p-1}
\end{equation}
or equivalently
\begin{equation}
    \nu_p\left[H(n)\right]\geq p\Bigl\lfloor\frac{n}{p}\Bigl\rfloor\left(\frac{n}{p-1}-\Bigl\lfloor\log_p(n)\Bigl\rfloor-1\right)-\frac{p}{2(p-1)}\left(\Big\lfloor\frac{n}{p}\Bigl\rfloor-1\right)\left(p\Big\lfloor\frac{n}{p}\Bigl\rfloor-2\right),
\end{equation}
which can be further simplified into
\begin{equation}
    \nu_p\left[H(n)\right]\geq \frac{p}{2(p-1)}\left\{-2+\Bigl\lfloor\frac{n}{p}\Bigl\rfloor\left[4+2n+p\left(1+\Bigl\lfloor\frac{n}{p}\Bigl\rfloor\right)\right]\right\}-p\Bigl\lfloor\frac{n}{p}\Bigl\rfloor\Bigl\lfloor\log_p(n)\Bigl\rfloor,
\end{equation}
as well as
\begin{equation}
    \nu_p\left[H(n)\right]\leq p\Big\lfloor\frac{n}{p}\Bigl\rfloor\left(\frac{n-1}{p-1}\right)-p\sum_{k=1}^{\Big\lfloor\frac{n}{p}\Bigl\rfloor-1}\left(2+\Bigl\lfloor\log_p(k)\Bigl\rfloor-\frac{pk}{p-1}\right)
\end{equation}
or equivalently
\begin{equation}
    \nu_p\left[H(n)\right]\leq p\Big\lfloor\frac{n}{p}\Bigl\rfloor\left(\frac{n-1}{p-1}\right)-p\sum_{k=1}^{\Big\lfloor\frac{n}{p}\Bigl\rfloor-1}\Bigl\lfloor\log_p(k)\Bigl\rfloor+\frac{p}{2(p-1)}\left(\Big\lfloor\frac{n}{p}\Bigl\rfloor-1\right)\left(4-4p+p\Big\lfloor\frac{n}{p}\Bigl\rfloor\right),
\end{equation}
which can be simplified into
\begin{equation}
    \nu_p\left[H(n)\right]\leq \frac{p}{2(p-1)}\left\{4(p-1)+\Big\lfloor\frac{n}{p}\Bigl\rfloor\left[2+2n+p\left(\Big\lfloor\frac{n}{p}\Bigl\rfloor-5\right)\right]\right\}-p\sum_{k=1}^{\Big\lfloor\frac{n}{p}\Bigl\rfloor-1}\Bigl\lfloor\log_p(k)\Bigl\rfloor.
\end{equation}

\subsection{Superfactorial}

The superfactorial is defined as
\begin{equation}
    \mathrm{sf}(n)=\prod_{k=1}^nk!.
\end{equation}
The valuation of the superfactorial is thus
\begin{equation}
    \nu_p[\mathrm{sf}(n)]=\sum_{k=1}^n\nu_p(k!)
\end{equation}
and using Eq. (\ref{majo1}), we have simply
\begin{equation}
    \nu_p[\mathrm{sf}(n)]\leq\sum_{k=1}^n\frac{k-1}{p-1}=\frac{n(n-1)}{2(p-1)}.
\end{equation}

\section{The Berezin function and further definitions}\label{sec3}

Berezin discussed the mathematical implications of the function defined by \cite{Berezin1987,Pickover1990,Pickover2001}:
\begin{equation}
    n\$=n!^{n!^{n!^{\cdot^{\cdot^{\cdot}}}}}.
\end{equation}
The term $n!$ is repeated $n!$ times. Then basically there is a ``tower'' of $n$ a's. In other words, there are $n!$ repetitions of $n!$ in the right-hand side. For example,
\begin{equation}
3\$=(3!)^{3!^{3!^{3!^{3!^{3!}}}}}.
\end{equation}
The function grows very rapidly and as $n$ increases, It reads $n\$=n!\uparrow\uparrow n!$ in Knuth's up-arrow\footnote{In the Knuth notation, $a\uparrow\uparrow n$ means raise $a$ to itself $n-1$ times. For example,
\begin{equation}
  a\uparrow\uparrow 4=a^{a^{a^a}}.  
\end{equation}
Such an operation is sometimes written $^na$ and named ``tetration'' \cite{Goodstein1947}.}.
It is worth mentioning that Aebi and Cairns proposed generalizations of Wilson's theorem for double-, hyper-, sub- and superfactorials \cite{Aebi2015}\footnote{For instance, if $p$ is an odd prime, then modulo $p$: $\mathrm{sf}(p-1)\equiv(-1)^{\frac{p-1}{2}}\,H(p-1)$.} In addition, the double hyperfactorial (OEIS sequence A002109 \cite{OEIS}) is defined as

\begin{equation}
    H_2(n)=\left\{
    \begin{array}{l}
    n^n\cdot(n-2)^{n-2}\cdots 5^5\cdot 3^3\cdot 1^1, \;\;\;\; \mathrm{for} \;\;\;\; n>0 \;\;\;\; \mathrm{and} \;\;\;\; n \;\;\;\;\mathrm{odd},\\
    n^n\cdot(n-2)^{n-2}\cdots 6^6\cdot 4^4\cdot 2^2, \;\;\;\; \mathrm{for} \;\;\;\; n>0 \;\;\;\; \mathrm{and} \;\;\;\; n \;\;\;\; \mathrm{even},\\
    0, \;\;\;\; \mathrm{for} \;\;\;\; n=0.
    \end{array}
    \right.
\end{equation}
which can be summarized as
\begin{equation}
    H_2(n)=\prod_{k=0}^{\lfloor\frac{n-1}{2}\rfloor}(n-2k)^{n-2k}
\end{equation}
and satisfies the recurrence
\begin{equation}
    H_2(n)=n^n\,H_2(n-2)\;\;\;\;\mathrm{with}\;\;\;\;H_2(0)=H_2(1)=1.
\end{equation}
The double hyperfactorial satisfies the properties
\begin{equation}
H_2(n)=\frac{1}{H_2(n-1)}\sqrt{\frac{H_2(2n)}{2^{n(n+1)}}}
\end{equation}
and
\begin{equation}
    H_2(n)H_2(n-1)=H(n).
\end{equation}
Other quantities can also be considered, such as 
\begin{equation}
    F_1(n)=\prod_{k=1}^nk^{k!},
\end{equation}
or
\begin{equation}
    F_2(n)=\prod_{k=1}^nk!^{k},
\end{equation}
or again
\begin{equation}
    F_3(n)=\prod_{k=1}^nk!^{k!}.
\end{equation}
Note that the function defined by $U(n)=(n!)^{n!}$ is sometimes referred to as the ultrafactorial (OEIS sequence A046882 \cite{OEIS}). Its values for $n$=0, 1, 2, 3 4 and 5 are respectively 1, 1, 4, 46656 and 1333735776850284124449081472843776.

We have $\forall \,n\geq 1$:
\begin{equation}
n!=\prod_{\substack{i+j=n+1\\i,j>1}}\sqrt{i\,j}
\end{equation}
and $\forall\, i,j\geq 1$:
\begin{equation}
i+j-1\leq i\,j\leq \left(\frac{i+j}{2}\right)^2
\end{equation}
which implies
\begin{equation}
n^{n/2}\leq n!\leq\left(\frac{n+1}{2}\right)^n.
\end{equation}
Thus, we have some bounds for the function $F_3(n)$:

\begin{equation}
\prod_{k=1}^nk^{k^{k/2}}\leq F_1(n)\leq\prod_{k=1}^nk^{\left(\frac{k+1}{2}\right)^k}
\end{equation}
and
\begin{equation}
\prod_{k=1}^nk^{k^2/2}\leq F_2(n)\leq\prod_{k=1}^n\left(\frac{k+1}{2}\right)^{k^2}
\end{equation}
as well as
\begin{equation}
\prod_{k=1}^nk^{\frac{1}{2}k^{1+\frac{k}{2}}}\leq F_3(n)\leq\prod_{k=1}^n\left(\frac{k+1}{2}\right)^{k\left(\frac{k+1}{2}\right)^k}.
\end{equation}
Concerning the $p$-adic valuations, one has
\begin{equation}
\nu_p\left[F_1(n)\right]=\sum_{k=1}^nk!\,\nu_p(k)
\end{equation}
but since
\begin{equation}
k\geq p^{\nu_p(k)}>2^{\nu_p(k)}>\nu_p(k), 
\end{equation}
one has
\begin{equation}
\nu_p\left[F_1(n)\right]\leq\sum_{k=1}^nk!\times k
\end{equation}
and since
\begin{equation}\label{kkfact}
\sum_{k=1}^nk!\times k=\sum_{k=1}^nk!\times (k+1)-\sum_{k=1}^nk!=(n+1)!-1
\end{equation}
one can write
\begin{equation}
\nu_p\left[F_1(n)\right]\leq (n+1)!-1.
\end{equation}
Of course, one has, since $k\geq p^{\nu_p(k)}$, the better upper bound
\begin{equation}
    \nu_p(k)\leq\log_p(k), 
\end{equation}
and then
\begin{equation}
\nu_p\left[F_1(n)\right]\leq\frac{1}{\log p}\sum_{k=1}^nk!\times \log k,
\end{equation}
but the sum can not be easily simplified. Similarly
\begin{equation}
\nu_p\left[F_2(n)\right]=\sum_{k=1}^nk\,\nu_p(k!)
\end{equation}
and thus, again using inequality (\ref{majo1}), one gets
\begin{equation}
\nu_p\left[F_2(n)\right]\leq\sum_{k=1}^n\frac{k(k-1)}{(p-1)}=\frac{n(n^2-1)}{3(p-1)}.
\end{equation}
For the lower bound, we have
\begin{equation}
    \nu_p\left[F_2(n)\right]\geq \sum_{k=1}^n\left(\frac{k}{p-1}-\Bigl\lfloor\log_p k\Bigl\rfloor-1\right)k,
\end{equation}
which can be simplified into
\begin{equation}
    \nu_p\left[F_2(n)\right]\geq \frac{n(n+1)(2n+4-3p)}{6(p-1)}-\sum_{k=1}^nk\Bigl\lfloor\log_p k\Bigl\rfloor.
\end{equation}

Similarly, for $F_3(n)$, one gets
\begin{equation}
\nu_p\left[F_3(n)\right]=\sum_{k=1}^nk!\,\nu_p(k!)
\end{equation}
and thus, again using inequality (\ref{majo1}), one gets
\begin{equation}
\nu_p\left[F_3(n)\right]\leq\sum_{k=1}^n\frac{k!\times k}{(p-1)}
\end{equation}
and using relation (\ref{kkfact}), one gets
\begin{equation}
\nu_p\left[F_3(n)\right]\leq\frac{\left[(n+1)!-1\right]}{(p-1)}.
\end{equation}
Taking into account the better upper bound for $\nu_p(k!)$ (see Eq. (\ref{mainineq})) yields
\begin{equation}
\nu_p\left[F_3(n)\right]\leq\sum_{k=1}^n\frac{k!\times (k-1)}{(p-1)}
\end{equation}
but the result is more complicated since it reads
\begin{equation}\label{f33}
\nu_p\left[F_3(n)\right]\leq!(-1) + \Gamma[2 + n] \frac{1+(-1)^n !(-2 - n))}{p-1},
\end{equation}
where $!n$ (also denoted $d(n)$), sometimes refereed to as ``subfactorial'', gives the number of permutations of $n$ objects that leave no object fixed (i.e., ``derangements''). In Eq. (\ref{f33}), the analytic continuation of the analytic continuation $!$ to negative values is used. The first values, starting from $n=1$ are 0, 1, 2, 9, 44, 265, 1854, 14833, 133496, 1334961, $\cdots$. such a number is equal to (see for instance Refs. \cite{Hassani2003,Stanley2012}):
\begin{equation}
!n=n!\sum_{k=0}^n\frac{(-1)^k}{k!}=\sum_{k=0}^nk!(-1)^{n-k}\binom{n}{k}=\sum_{k=0}^n\frac{n!(-1)^{n-k}}{(n-k)!}=\frac{\Gamma(n+1,-1)}{e}
\end{equation}
and $\Gamma(a,z)$ is the incomplete Gamma function. The asymptotic form is \cite{Hassani2020}:
\begin{equation}
!n={\frac {n!}{e}}+\sum _{k=1}^{m}\left(-1\right)^{n+k-1}{\frac {B_{k}}{n^{k}}}+O\left({\frac {1}{n^{m+1}}}\right),
\end{equation}
where $m$ is any fixed positive integer, and $B_{k}$ denotes the $k$-th Bell number. The following recurrence also holds:
\begin{equation}
!n={\begin{cases}1&{\text{if }}n=0,\\n\cdot \left(!(n-1)\right)+(-1)^{n}&{\text{if }}n>0.\end{cases}}.
\end{equation}
For the lower bound, one has
\begin{equation}
    \nu_p\left[F_3(n)\right]\geq \sum_{k=1}^n\left(\frac{k}{p-1}-\Bigl\lfloor\log_p k\Bigl\rfloor-1\right)k!,
\end{equation}
which can be simplified into
\begin{align}
    \nu_p\left[F_2(n)\right]\geq &\frac{(p-2)+(p-1)(!(-1))+(n+1)!\left[1+(-1)^n)(p-1)(!(-2-n))\right]}{p-1}\\
    &-\sum_{k=1}^nk!\Bigl\lfloor\log_p k\Bigl\rfloor,
\end{align}
where the analytic continuation of $!n$ to negative values is used as well.

\section{Conclusion}

We obtained simple bounds for the $p$-adic valuation $\nu_p$ of the factorial $n!$, the hyperfactorial $H(n)$, the superfactorial $\mathrm{sf}(n)$, the Berezin-Pickover function and other similar quantities (involving factorials, powers of factorials or iterated factorials) likely to yield very large numbers. Such functions,  are important for cryptography and data protection. The derivations are based on the Legendre-de Polignac formula. I addition, lower and upper bounds for the $p$-adic valuation of Stirling numbers of the first kind and Catalan numbers are given. In the future, we plan to derive bounds for the primorial, hyperprimorial and superprimorial functions \cite{Rosser1962,Dubner1987,Mezo2013}, and to use the expression of the $p$-adic valuation in terms of binomial coefficient to obtain new identities.

\appendix

\section{Appendix: Other interesting results concerning the $p$-adic valuation}\label{appA}

In this appendix, we provide expressions of the $p$-adic valuation which can be useful in order to derive new relations or bounds. De Castro recently found that, for a positive integer $n$ and a prime number $p$, the $p$-adic valuation of $n$ can be expressed as
\begin{equation}
    \nu_p(n)=p\sum_{j=1}^{\lfloor\log_p n\rfloor}\left\{\binom{n}{p^j}\frac{p^{j-1}}{n}\right\},
\end{equation}
which can alternatively be expressed with the use of complex numbers \cite{Castro2022}:
\begin{equation}
    \nu_p(n)=\frac{i\,e^{-i\pi/p}}{2\sin(\pi/p)}\left(\lfloor\log_p n\rfloor-\sum_{j=1}^{\lfloor\log_p n\rfloor}e^{\frac{2i\pi}{p}\binom{n-1}{p^j-1}}\right).
\end{equation}
One has
\begin{equation}
\binom{n}{k}=\frac{n(n-1)\cdots (n-k+1)}{k!}\leq\frac{n^k}{k!}-1,
\end{equation}
and
\begin{equation}
e^k=\sum_{l=0}^{\infty}\frac{k^l}{l!}>\frac{k^k}{k!}
\end{equation}
yielding $1/k!<(e/k)^k$ and thus
\begin{equation}
    \binom{n}{k}\leq\left(\frac{en}{k}\right).
\end{equation}
Considering an integer $m$ such that $0<m<k\leq n$, one has
\begin{equation}
    \frac{m}{n}\leq\frac{m}{k}\Rightarrow \frac{k-m}{k}\leq\frac{n-m}{n}\Rightarrow \frac{n}{k}\leq\frac{n-m}{k-m}
\end{equation}
and thus 
\begin{equation}
    \left(\frac{n}{k}\right)^k=\frac{n}{k}\cdots \frac{n}{k}\leq \frac{n}{k}\cdot\frac{(n-1)}{(k-1)}\cdots\frac{(n-k+1)}{1}
\end{equation}
so that we have the well-known inequalities
\begin{equation}\label{wk}
    \left(\frac{n}{k}\right)^k\leq\binom{n}{k}\leq\left(\frac{en}{k}\right)^k.
\end{equation}
We can therefore write
\begin{equation}
    \sum_{j=1}^{\lfloor\log_p n\rfloor}\left(\frac{n}{p^j}\right)^{p^j}\frac{p^j}{n}\leq \nu_p(n)\leq\sum_{j=1}^{\lfloor\log_p n\rfloor}\left(\frac{en}{p^j}\right)^{p^j}\frac{p^j}{n}.
\end{equation}
Such an inequality can also be of interest in order to obtain new bounds for $F_1(n)$ (see section \ref{sec3}).

Agievich recently published the better upper bound \cite{Agievich2022}:
\begin{equation}\label{agie}
    \binom{n}{k}\leq\frac{2^n}{\sqrt{\pi n/2}}\,\exp\left[-\frac{2}{n}\left(k-\frac{n}{2}\right)^2+\frac{23}{18n}\right],
\end{equation}
which thus gives
\begin{equation}
    \nu_p(n)\leq\sum_{j=1}^{\lfloor\log_p n\rfloor}\frac{2^{p^j}}{\sqrt{\frac{\pi p^j}{2}}}\frac{p^j}{n}\,\exp\left[-\frac{2}{n}\left(p^j-\frac{n}{2}\right)^2+\frac{23}{18n}\right],
\end{equation}
also relevant for inequalities involving $F_1(n)$.

Komatsu and Young investigated the $p$-adic valuation of Stirling numbers of the first kind \cite{Komatsu2017}. The latter are defined through the generating function
\begin{equation}
    x(x-1)\cdots(x-n+1)=\sum_{k=0}^ns(n,k)\,x^k.
\end{equation}
Let $k$ be a positive integer and $p$ a prime. If $n$ is of the form $n=kp^r+m$, where $0\leq m<p^r$, then
\begin{equation}
    \nu_p\left[s(n+1,k+1)\right]=\nu_p(n!)-\nu_p(k!)-kr.
\end{equation}
Equivalently, one can write
\begin{equation}\label{equistir}
    \nu_p\left[s(n+1,k+1)\right]=k\left(\frac{p^r-1}{p-1}-r\right)+\nu_p(m!).
\end{equation}
Defining, for $k$ a non-negative integer and $p$ a prime, the set $A_{k,p}=\left\{kp^r+m: r\geq 0, 0\leq m\leq p^r\right\}$, for all $n\in A_{k,p}$ with $n>k$ we have
\begin{equation}
    \frac{k}{n-k}\Bigl\lfloor\log_p\left(\frac{n}{k}\right)\Bigl\rfloor\leq\frac{1}{p-1}-\frac{\nu_p\left[s(n+1,k+1)\right]}{n-k}\leq\frac{k+1}{n-k}\Bigl\lfloor\log_p\left(\frac{n}{k}\right)\Bigl\rfloor.
\end{equation}
In addition, let a positive integer $n\geq k$ be given and define $r$ so that $kp^r\leq n\leq kp^{r+1}$. The authors found that
\begin{equation}
    \nu_p\left[s(n+1,k+1)\right]\geq\nu_p(n!)-\nu_p\left(\Bigl\lfloor\frac{n}{p^r}\Bigl\rfloor!\right)-kr.
\end{equation}

Using Eq. (\ref{equistir}) bounds from Eq. (\ref{mainineq}), one obtains
\begin{equation}\label{nara}
    k\left(\frac{p^r-1}{p-1}-r\right)+\frac{m}{p-1}-\Bigl\lfloor\log_p(m)\Bigl\rfloor-1\leq\nu_p\left[s(n+1,k+1)\right]\leq k\left(\frac{p^r-1}{p-1}-r\right)+\frac{m-1}{p-1}.
\end{equation}

Guadalupe characterized the 3-adic valuation of the Narayana sequence defined as $a_0=1$, $a_1=a_2=1$ and $a_n=a_{n-1}+a_{n-3}$ for $n\geq 3$ and used this to determine all Narayana numbers that are factorials \cite{Guadalupe2021}. He obtained in particular, using inequality (\ref{nara}), that for $m\geq 6$ and for some $\delta\in\left\{0,1,3,8\right\}$:
\begin{equation}
    \frac{m}{2}-\Big\lfloor\frac{\log m}{\log 3}\Big\rfloor-1\leq\nu_3(a_n)\leq\nu_3\left[n(n+1)(n+3)(n+8)\right]+6\leq 4\nu_3(n+\delta)+6.
\end{equation}
Recently also, Moll studied the 2-adic valuation of Catalan numbers:
\begin{equation}
    C_n=\frac{1}{n+1}\binom{2n}{n}.
\end{equation}
One has
\begin{equation}
\nu_2(C_n)=\nu_2\left[(2n)!\right]-2\,\nu_2(n!)-\nu_2(n+1).
\end{equation}
As noticed by Moll \cite{Moll2020}, since
\begin{equation}
\nu_2(n!)=\sum_{k=1}^{\infty}\Bigl\lfloor\frac{n}{2^k}\Bigl\rfloor=n-s_2(n),
\end{equation}
where $s_2(n)$ is the sum of the binary digits of $n$, one has
\begin{equation}
    \nu_2(C_n)=s_2(n)-\nu_2(n+1)=n-\nu_2\left[(n+1)!\right]=s_2(n+1)-1.
\end{equation}
and Moll obtained that $C_n$ is odd precisely when $n$ has only 1's in its binary expression. 

Using
\begin{equation}
\nu_p(C_n)=\nu_p\left[(2n)!\right]-2\,\nu_p(n!)-\nu_p(n+1)
\end{equation}
and applying (\ref{mainineq}), one gets
\begin{equation}
        \frac{2n}{p-1}-\Bigl\lfloor\log_p(2n)\Bigl\rfloor-1\leq\nu_p[(2n)!]\leq\frac{2n-1}{p-1}.
\end{equation}
as well as
\begin{equation}
    -\frac{n-1}{p-1}\leq-\nu_p[(n)!]\leq-\frac{n}{p-1}+\Bigl\lfloor\log_p(n)\Bigl\rfloor+1.
\end{equation}
and since one has
\begin{equation}
    -\log_p(n+1)\leq-\nu_p(n+1)\leq 0,
\end{equation}
one finally obtaind
\begin{equation}
    \frac{n+1}{p-1}-\Bigl\lfloor\log_p(2n)\Bigl\rfloor-1-\log_p (n+1)\leq\nu_p(C_n)\leq\frac{n-1}{p-1}+1+\Bigl\lfloor\log_p(n)\Bigl\rfloor.
\end{equation}
More generally, the $p$-adic valuation of $k$-central binomial coefficients was discussed by Straub \emph{et al.} \cite{Straub2009}.

\end{document}